\documentclass[11pt]{amsart}

\newcommand{\ti}{\textit}
\newcommand{\MB}{\mathcal B}
\newcommand{\MD}{\mathcal D}
\newcommand{\MY}{\mathcal Y}
\newcommand{\MX}{\mathscr X}
\newcommand{\MH}{\mathcal H}
\newcommand{\MZ}{\mathcal Z}
\newcommand{\MI}{\mathcal I}
\newcommand{\Ga}{\alpha}

\newcommand{\Ge}{\varepsilon}
\newcommand{\BZ}{\mathbb Z}
\newcommand{\BE}{\mathbb E}

\newcommand{\BN}{\mathbb N}

\newcommand{\bx}{\mathbf x}

\newcommand{\BX}{\mathbf X}
\newcommand{\BY}{\mathbf Y}

\newtheorem{theorem}{Theorem}[section]
\newtheorem{lemma}[theorem]{Lemma}
\newtheorem{proposition}[theorem]{Proposition}
\newtheorem{prop}[theorem]{Proposition}
\newtheorem{corollary}[theorem]{Corollary}

\newcommand{\varnorm}[1]{\lvert\!|\!| #1|\!|\!\rvert}

\theoremstyle{definition}

\newtheorem{remark}{Remark}

\newtheorem*{observation}{Observation}

\usepackage{amsmath}
\usepackage{enumerate}
\usepackage{amssymb}
\usepackage{mathrsfs}
\begin{document}

\title{Cubic ergodic averages for actions of amenable groups.}
\begin{abstract} We unify and extend some previous results about multiparameter configurations in sets of positive density in products of amenable groups.  We show that if $G$ is an amenable group, and $E\subset G^d$ has positive density with respect to some F{\o}lner sequence $\{\Phi_N\}_{N\in \BN}$ in $G^d,$ then $E$ contains many ``cubic" configurations.  This generalizes the author's result for the case $d=2,$ and Q. Chu's for the case $G=\BZ.$  \end{abstract}

\author{John T. Griesmer}

\date{\today}


\address{Department of Mathematics, The Ohio State University, 231 W. 18th Ave., Columbus, OH 43212}
\email{griesmer@math.osu.edu}

\maketitle

\section{Introduction}

Khintchine \cite{K} proved the following strengthening of Poincar\'e's recurrence theorem.
\begin{theorem}{\rm (Khintchine's Recurrence Theorem.)} Let $(X,\MX,\mu)$ be a probability space, and let $T:X\to X$ be an invertible transformation which preserves $\mu.$  For all $B\in \MX$ and all $c>0,$ the set $\{n\in \BZ: \mu(B\cap T^{-n}B)>\mu(B)^2-c\}$ has bounded gaps.
\end{theorem}
One can deduce from this the fact that if $A\subset \BZ$ has positive density, that is if $d^*(A):=\limsup_{N-M\to \infty} \frac{|A\cap [N,M]|}{M-N+1}>0,$ then for all $c>0,$ then the set $\{n: d^*(A\cap (A-n))> d^*(A)^2-c\}$ has bounded gaps.

Bergelson established the following natural two-dimensional generalization of Khintchine's recurrence theorem in \cite{B}.

\begin{theorem}\label{bergelson}  Let $(X,\MX,\mu)$ be a probability space, and let $T:X\to X$ be an invertible transformation which preserves $\mu.$  For all $B\in \MX$ and all $c>0,$ the set 
$$
\{(n,m): \mu(B\cap T^{-n}B\cap T^{-m}(B\cap T^{-n}B))>\mu(B)^4-c\}
$$ meets every large enough square $[N_1,N_1+M]\times [N_2,N_2+M]$ in $\BZ^2.$
\end{theorem}

Bergelson applied Theorem \ref{bergelson}, together with the Furstenberg correspondence principle, to deduce
\begin{corollary}\label{bergelson2}  If $A\subset \BZ$ with 
$
d^*(A)>0,$ then 
$$
\{(n,m):d^*(A\cap  (A-n)\cap (A-m)\cap (A-n-m))>\mu(A)^4-c\}
$$ meets every large enough square $[N_1,N_1+M]\times [N_2,N_2+M]$ in $\BZ^2.$ 
\end{corollary}
Roughly speaking, this says that if $A\subseteq \BZ$  with $d^*(A)=\delta>0,$ then there are many $(n,m)\in \BZ^2$ for which the quadruples $\{a,a+n,a+m,a+n+m\}$ appear in $A$  (almost) as often as one would expect in a set $A$ generated randomly by selecting each $a\in \BZ$ independently with probability $\delta.$ 

Higher-dimensional generalizations of Theorem \ref{bergelson} and Corollary \ref{bergelson2} were obtained by Host and Kra in \cite{HK}.

One may wonder if results similar to Corollary \ref{bergelson2} hold for subsets of $\BZ^2$.  Defining the \ti{upper Banach density} of a subset $A\subset \BZ^2$ as 
$$
d^*(A):=\limsup_{\min(M_1-N_1,M_2-N_2)\to \infty} \frac{A\cap([N_1,M_1]\times [N_2,M_2])}{(M_1-N_1+1)(M_2-N_2+1)},
$$ one may ask if $d^*(A\cap (A-(n,0))\cap (A-(0,m))\cap (A-(n,m)))> d^*(A)^4-c$ for many $n,m$.  More generally, one may ask if similar results hold for subsets of an arbitrary group $G$.  

Of course, the question implies that there is a notion of density on $G$ similar to that of $d^*$ on $\BZ$ and $\BZ^2.$  As a substitute for sequences of intervals $[N,M],$ we use F{\o}lner sequences, which we define now.

If $G$ is a countable discrete group, a sequence $\{\Phi_N\}_{N\in \BN}$ of finite subsets of $G$ is called a \ti{left} (resp.\ \ti{right}) \ti{F{\o}lner sequence} if for all $g\in G,$ 
$$
\lim_{N\to \infty} \frac{|\Phi_N\cap g \Phi_N|}{|\Phi_N|}=1 \text{\ (resp.\ }\lim_{N\to \infty} \frac{|\Phi_N\cap  \Phi_Ng|}{|\Phi_N|}=1).
$$  A F{\o}lner sequence is called \ti{two-sided} if it is both a left- and a right F{\o}lner sequence.  If a discrete group $G$ has a F{\o}lner sequence, then $G$ is called \ti{amenable}.  We will usually denote F{\o}lner sequences $\{\Phi_N\}_{N\in \BN}$ without subscripts, as in ``let $\Phi$ be a F{\o}lner sequence."

If $\Phi$ is a left F{\o}lner sequence in a group $G,$ we can define the upper density (with respect to $\Phi$) of a set $A\subset G$ with respect to $\Phi$ by $d_{\Phi}(A)=\limsup_{N\to \infty} \frac{|A\cap \Phi_N|}{|\Phi_N|}.$ It is easy to check that $d_{\Phi}(A)=d_{\Phi}(gA)$ for all $g\in G$ and $A\subseteq G.$  One can define density with respect to a right F{\o}lner sequence in the same way, so that $d_\Phi$ is invariant under right multiplication if $\Phi$ is a right F{\o}lner sequence.

We also need a notion that generalizes the idea of a subset $R$ of $\BZ^2$ that meets every large enough square $[N_1,N_1+M]\times [N_2,N_2+M].$  We call a subset $S\subset G$ \textit{left- (resp.\ right-) syndetic} if there exists finitely many $g_1,\dots, g_k\in G$ such that $G=\bigcup_{i=1}^k g_iS$ (resp.\ $G=\bigcup_{i=1}^k Sg_i$).  One can easily check that $A\subset \BZ$ is syndetic if and only if $A$ meets every long enough interval, and $A\subset \BZ^2$ is syndetic if and only if $A$ meets every large enough square in $\BZ^2.$

The following generalization of Corollary \ref{bergelson2} appears in \cite{Gr}\footnote{The exponents $-1$ do not appear in \cite{Gr}, but they can be removed or added at will via a change of variables.}.  Here $E(a,b)$ means $\{(xa,yb):(x,y)\in E\}.$  
\begin{theorem}\label{Grcor}{\rm(\cite{Gr}, Corollary 1.6)} Let $G$ be a countable amenable group with identity $e.$  Let $E\subseteq G\times G,$ and let $\Phi$ be a right F{\o}lner sequence in $G\times G$ with $\limsup_{N\to \infty} \frac{|E\cap \Phi_N|}{|\Phi_N|}=\delta>0.$  Then for all $c>0,$ the set
\begin{align*}
\left\{(g,h): \limsup_{N\to \infty} \frac{|E\cap E(g^{-1},e)\cap E(e,h^{-1})\cap E(g^{-1},h^{-1})\cap\Phi_N|}{|\Phi_N|}>\delta^4-c\right\}
\end{align*}
is both left- and right syndetic in $G\times G.$
\end{theorem}

This is deduced from the following theorem and corollary in \cite{Gr}.  (We say two actions $T,S$ of $G$ commute if $S_gT_h=T_hS_g$ for all $g,h\in G$.)

\begin{theorem}\label{Grthm} {\rm(\cite{Gr}, Theorem 1.4)}  Let $G$ be a countable amenable group, let $\mathbf X=(X,\MX, \mu)$ be a probability space with  probability measure $\mu,$ and commuting measure preserving actions $T,S$ of $G$ on $X.$  Then

(1) For all $f_1,f_2,f_3\in L^\infty (\mu),$ and all two-sided F{\o}lner sequences $\Phi,\Psi$ in $G,$ the limit
\begin{align*}
L=\lim_{N\to \infty} \frac{1}{|\Phi_N||\Psi_N|} \sum_{(g,h)\in \Phi_N \times \Psi_N} f_1(T_g x) f_2(S_hx) f_3(T_gS_hx)
\end{align*}
exists in $L^2(\mu).$  

(2) $L$ is equal to a constant $\mu$-almost everywhere for all $f_1,f_2,f_3\in L^\infty(\mu)$ if and only if $T\times T$ and $S\times S$ are ergodic. 
\end{theorem}
It is also shown in \cite{Gr} that the limit $L$ above is independent of the choice of F{\o}lner sequences $\Phi,\Psi.$

\begin{corollary}\label{Grbound}{\rm(\cite{Gr}, Corollary 1.5)} Let $(X,\mathcal B,\mu),T$ and $S$ be as in Theorem \ref{Grthm}.  Suppose $f\in L^\infty(\mu)$ is a nonnegative function.  Then  for all two-sided F{\o}lner sequences $\Phi,\Psi$ 
\begin{align}\label{large1}
\lim_{n\to \infty} \frac{1}{|\Phi_N||\Psi_N|} \sum_{(g,h)\in \Phi_N\times \Psi_N} \int fT_gfS_hfT_gS_hf\, d\mu\geq \left(\int f\, d\mu \right)^4.
\end{align}
\end{corollary}

In this paper we generalize Theorem \ref{Grthm} and Theorem \ref{Grcor} to higher-dimensional cases.  Specifically, we establish the following theorem and corollary.
\begin{theorem}\label{main} Let $G$ be a countable amenable group, let $(X,\mathcal B, \mu)$ be a probability space with  probability measure $\mu,$ let $d\in \BN,$ and let $T^{(1)}, \dots, T^{(d)}$ be commuting actions of $G$ on $X$ which preserve $\mu.$  Let $\Phi^{(1)},\dots, \Phi^{(d)}$ be F{\o}lner sequences in $G,$ and let $F_N$ denote $\Phi^{(1)}\times \cdots \times \Phi^{(d)}.$  Let $f_{\Ge}, \Ge\in \{0,1\}^d$ be $2^d$ bounded, $\MX$-measurable functions on $X.$   
\begin{enumerate} \item[(1)]
The limit
$$
L=\lim_{N\to \infty} \frac{1}{|F_N|} \sum_{g\in F_N} \prod_{\Ge\in \{0,1\}^d} f_\Ge \circ R_g^{\Ge}
$$
exists in $L^2(\mu),$ where $R_g^{\Ge}:= \prod_{i=1}^d(T^{(i)}_{g_i})^{\Ge_i}$ for $g=(g_1,\dots, g_d), \Ge=(\Ge_1,\dots, \Ge_d).$

\item[(2)] Furthermore, if each $f_\Ge=f$ is a nonnegative function, then
$$
\int L\, d\mu \geq \left(\int f\, d\mu \right)^{2^d}.
$$
\end{enumerate}

\end{theorem}

\begin{corollary}\label{maincor}
Let $G$ be a countable amenable group, let $G^d$ denote the $d^\text{th}$ cartesian power of $G,$ and let $E\subset G^d.$  For $g=(g_1,\dots, g_d), h=(h_1,\dots, h_d)\in G^{d},$ let $
h.g:=\{(h_1^{-\Ge_1}g_1,\cdots, h_d^{-\Ge_d}g_d): \Ge\in \{0,1\}^d\}.$
If there is a left F{\o}lner sequence $\Psi$ in $G^d$ with $d_\Psi(A)>0,$ then for all $c>0,$ the set
$$
\{h\in G^d: d_\Psi(\{g:h.g\subset E\})>d_\Psi(E)^{2^d}-c\}
$$
is both left- and right- syndetic in $G^d.$
\end{corollary}
\begin{remark}
One may make the analogous conclusion about right F{\o}lner sequences if one replaces $h.g$ with $g_*h:=\{(g_1h_1^{\Ge_1},g_2h_2^{\Ge_2},\cdots, g_dh_d^{\Ge_d}): \Ge\in \{0,1\}^d\}.$

Note that $h.g\subseteq E$ if and only if $g\in \bigcap_{\Ge\in \{0,1\}^d} (h_1^{\Ge_1},\dots, h_d^{\Ge_1})E,$ so the conclusion is equivalent to concluding that 
$$
\{h\in G^d: d_\Psi\left(\bigcap_{\Ge\in \{0,1\}^d} (h_1^{\Ge_1},\dots, h_d^{\Ge_d}) E\right)>\delta^{2^d}-c\}
$$ is both left- and right syndetic in $G^d.$
\end{remark}
Theorem \ref{main} and Corollary \ref{maincor} were established for the case $G=\BZ$ in \cite{C}, using results from \cite{H}.  Theorem \ref{main} for the case $G=\BZ$ is Theorem 1.1 of \cite{C}. Although Theorem \ref{main} was established for the case $G=\BZ$ in \cite{A}, the methods of \cite{H} and \cite{C} provide additional insight.

\subsection{Acknowledgments.}  The author would like to thank Vitaly Bergelson for helpful comments and advice.

\section{Preparation the Proof of Theorem \ref{main}.}

The papers \cite{C} and \cite{H} together constitute a self-contained proof of Theorem \ref{main} in the case $G=\BZ.$  A self-contained proof of Theorem \ref{main} for an arbitrary amenable group $G$ would be virtually identical to the papers \cite{C} and \cite{H}, except that the proofs of Lemma 2 in \cite{H} and Theorem 4 of \cite{H} must be changed, and some routine facts about F{\o}lner sequences must be applied.

Before we present those changes, we recall some background.

\subsection{Background.}

By a standard probability space, we mean a measure space $(X,\MX,\mu)$ which is measure-theoretically isomorphic to a measure space $(K,\MB, \nu)$ where $K$ is a compact metric space, $\MB$ is the $\sigma$-algebra of Borel subsets of $K,$ and $\nu$ is a regular Borel measure on $\MB$ with $\nu(K)=1.$  The relevant probability spaces for deducing Corollary \ref{maincor} have compact metric topologies on the underlying spaces, and the general case of Theorem \ref{main} can be deduced from the special case where $(X,\MX,\mu)$ is a standard probability space by an application of Theorem 2.15 of \cite{G}.

Let $(X,\MX,\mu)$ be a standard probability measure space, let $G$ be a countable group, and let $T$ be an action of $G$ on $X$ which preserves $\MX$ and $\mu.$  We say that $\BX=(X,\MX,\mu,T)$ is a measure preserving $G$-system (or $G$-system).  

If $k\in \BN$ , a \ti{$k$-fold self-joining} of $\BX$ is a $G$-system $(X^k, \MX^{\otimes k}, \lambda, T^{\otimes k}),$ where $X^k$ is the $k^\text{th}$ cartesian power of $X, \MX^{\otimes k}$ is the $k$-fold product $\sigma$-algebra, $T^{\otimes k}$ is the action given by $g\mapsto T^{\otimes k}_g,$ where $T^{\otimes k}$ is the $k$th cartesian power of $T,$ and $\lambda$ is a $T^{\otimes k}$-invariant probability measure on $(X^k,\MX^{\otimes k})$ satisfying the following condition: if $A\in \MX,$ then $\lambda(\prod_{i=1}^{j-1} X \times A\times \prod_{i=j+1}^k X )= \mu(A),$ for $1\leq j\leq k.$

If $\BY=(Y,\MY,\nu, S)$ is another $G$-system, we say that $\BY$ is a \ti{factor} of $\BX$ if there is a map $\pi:X\to Y$ with $\pi^{-1}(\MY)\subseteq \MX, \mu(\pi^{-1}(A))=\nu(A)$ for all $A\in \MY,$ and $\pi \circ T_g=S_g\circ \pi$ for all $g\in G.$

If $\MD\subset \MX$ is any countably generated, $T$-invariant sub $\sigma$-algebra, one can always realize $\MD$ as a $\sigma$-algebra of the form $\pi^{-1}(\MY),$ where $(Y,\MY, \nu, S)$ is a factor of $\MX$ and $(Y,\MY,\nu)$ is a standard probability space.  In this way, we have a correspondence between factors of $\BX$ and $T$-invariant sub $\sigma$-algebras of $\MX.$

If $\MD\subset \MX$ is a countably generated sub $\sigma$-algebra of $\MX,$ and $f\in L^2(\mu),$ we let $\BE_\mu(f|\MD)$ denote the conditional expectation of $f$ on $\MD,$ which agrees with the orthogonal projection of $f$ onto the closed subspace of $L^2(\mu)$ consisting of the $\MD$-measurable functions.  If there is no ambiguity , we write $\BE(f|\MD)$ for $\BE_\mu(f|\MD).$

If $\BY$ is a factor of $\BX,$ the \ti{relatively independent product} of $\BX$ with itself over $\BY$ is the $2$-fold self joining $(X\times X, \MX\otimes \MX, \mu \times_\MY \mu, T\times T),$ where $\mu\times_\MY\mu$ is given by $\int f\otimes g\, d\mu\times_\MY\mu= \int \BE(f|\MY) \BE(g|\MY)\, d\mu,$ for $f,g\in L^\infty(\mu).$

If $(X,\MX,\mu, T)$ is a $G$-system, we let $\MI_T$ denote the $\sigma$-algebra of $T$-invariant $\MX$-measurable sets.

\subsubsection{The mean ergodic theorem.}

The mean ergodic theorem for amenable groups follows.

\begin{theorem}\label{ergodictheorem}  Let $(X,\MX,\mu)$ be a probability space, let $G$ be an amenable group, $T:X\to X$ be a $\mu$-preserving action of $G$ on $X.$  Let $\Phi$ be a left F{\o}lner sequence in $G.$  Then for all $f\in L^2(\mu),$
$$
\lim_{N\to \infty} \frac{1}{|\Phi_N|} \sum_{g\in \Phi_N} f\circ T_g= Pf
$$
in $L^2(\mu),$ where $P$ is the orthogonal projection of $f$ onto the space of $T$-invariant functions. 
\end{theorem}
We will need a slight variation, which follows immediately from Theorem \ref{ergodictheorem}.

\begin{proposition}\label{variation}
Let $(X,\MX,\mu)$ be a probability space, let $G$ be an amenable group, $T:X\to X$ be a $\mu$-preserving action of $G$ on $X.$  Let $\Phi$ be a right F{\o}lner sequence in $G.$  Then for all $f\in L^2(\mu),$
$$
\lim_{N\to \infty} \frac{1}{|\Phi_N|^2} \sum_{g,h\in \Phi_N} f\circ T_{gh^{-1}}= Pf$$
in $L^2(\mu),$ where $P$ is the orthogonal projection of $f$ onto the space of $T$-invariant functions.
\end{proposition}
To prove this, write $Pf=\frac{1}{|\Phi_N|} \sum_{h\in \Phi_N} Pf\circ T_{h},$ so that $$
\frac{1}{|\Phi_N|^2} \sum_{g,h\in \Phi_N} f\circ T_{gh^{-1}}-Pf= \frac{1}{|\Phi_N|} \sum_{h\in \Phi_N} \frac{1}{|\Phi_N|} \sum_{g\in \Phi_N} (f\circ T_{g}-Pf)\circ T_{h^{-1}}.
$$
Taking norms, we have
\begin{align*}
\left\|\frac{1}{|\Phi_N|^2} \sum_{g,h\in \Phi_N} f\circ T_{gh^{-1}}-\right. & \left. Pf\right\|_{L^2(\mu)} \\ &\leq \frac{1}{|\Phi_N|} \sum_{h\in \Phi_N} \left\|\frac{1}{|\Phi_N|} \sum_{g\in \Phi_N} (f\circ T_{g}-Pf)\circ T_{h}\right\|\\
&= \left\|\frac{1}{|\Phi_N|} \sum_{g\in \Phi_N} (f\circ T_{g}-Pf)\right\|.
\end{align*}
Now apply Theorem \ref{ergodictheorem}.

\subsection{F{\o}lner sequences and syndeticity.}

To derive Corollary \ref{maincor} from Theorem \ref{main}, we need the following lemma, whose proof is so similar to the proof of Lemma 4.7 of \cite{Gr} that we omit it.

\begin{lemma}\label{tosynd} Let $G$ be an amenable group, $d\in \BN,$ and let $S\subset G^d.$  Then $S$ is both left and right syndetic if and only if for all two-sided F{\o}lner sequences $\Phi^{(i)},i=1,\dots, d$ in $G,$ there exists $N_1,\dots, N_d\in \BN$ such that $S\cap \Phi^{(1)}_{N_1}\times \cdots \times \Phi^{(d)}_{N_d}\neq \emptyset.$
\end{lemma}

\section{Proof of Theorem \ref{main}.}

We will present the constructions of \cite{H}, adapted to the setting of our Theorem \ref{main}, and then explain how these are used in \cite{C} to establish Theorem 1.1 of \cite{C}.

\subsection{The cubes joinings.}

Throughout this section, we fix a standard probability space $(X,\MX,\mu),$ an amenable group $G,$ and commuting, $\mu$-preserving actions of $G$ on $(X,\MX),$ denoted by $T^{(1)},\dots, T^{(d)}.$  Let $\MI_i$ denote the $\sigma$-algebra consisting of $T^{(i)}$-invariant sets.

Unless otherwise stated, for $k>0,$ we consider the set $\{0,1\}^{k}$ to have the reverse-lexicographic order, and we write $\mathbf 0$ for the element of $\{0,1\}^k$ all of whose coordinates are $0.$

For each sequence $P=(i_1,i_2,\dots, i_k),$ of distinct integers $i_j\leq d,$  we define a joining $\mu^P$ on $X^{2^k}$ as follows.  For $P=\emptyset, \mu^P=\mu.$  If $k>0, P=(i_1,\dots, i_k)$  and $\mu^{Q}$ is defined for all $Q$ of length $k-1,$  let $\MI_{Q,j}$ denote the $(T^{(j)})^{\otimes 2^{k-1}})$-invariant $\sigma$-algebra of the system $(X^{2^{k-1}},\mu^Q, (T^{(j)})^{\otimes 2^{k-1}}),$ let $Q=(i_1,\dots, i_{k-1}),$ and define $\mu^{P}$ on $X^{2^k}$ by
$$
\int f\otimes g\, d\mu^{P}= \int \BE(f|\MI_{Q,k})\BE(g|\MI_{Q,k})\, d\mu^Q,
$$
for $f,g\in L^\infty(X^{2^{k-1}}),$
so that $\mu^P$ is the relatively independent product of $\mu^Q$ with itself over $\MI_{Q,k}.$

By repeated application of the mean ergodic theorem, the joining $\mu^{P}$ can also be defined as follows.  For $f_{\Ge}\in L^\infty(X), \Ge\in {\{0,1\}^k},$ the integral $\int \bigotimes_{\Ge\in \{0,1\}^k} f_\Ge\, d\mu^{P}$ is equal to the iterated limit
$$
\lim_{N_1\to \infty} \dots \lim_{N_k\to \infty} \frac{1}{|\Phi_{N_1}|\cdots |\Phi_{N_k}|} \sum_{g\in \Phi_k} \cdots \sum_{g\in \Phi_i} \prod_{\Ge\in \{0,1\}^k} f_\Ge \circ \prod_{j=1}^k (T_{g_j}^{(j)})^{\Ge_j}\, d\mu.
$$

The next lemma appears as Theorem 4 in \cite{H} for the case $G=\mathbb Z,$ and is essential to the proof of Theorem 1.1 in \cite{C}.

\begin{lemma}  The joining $\mu^{(P)}$ for $P=(i_1,\dots, i_k)$ depends only on the indices $\{i_1,\dots, i_k\}$ and not on the order. 
\end{lemma}

\ti{Proof.}  We argue as in Proposition 3 of \cite{H}.  For the case $k=2,$ this was shown in \cite{Gr}.  We prove the general statement by induction on $k.$  The induction hypothesis is ``For $r<k$, the joining $\mu^{(i_1,\dots, i_r)}$ is independent of the order of the $(i_1,\dots i_r).$"  It follows that the joining $\mu^{(i_1,\dots, i_k)}$ is independent of the order of the first $r$ coordinates.  We then show that for $P=(i_1,\dots, i_k), Q=(i_2,i_1\dots, i_{k-1},i_k), \mu^P=\mu^Q.$

Write $U=T^{(i_{k-1})},V=T^{(i_k)}$ and let $U',V'$ denote the product actions $ (U)^{\otimes 2^{k-2}}, (V)^{\otimes 2^{k-2}},$ respectively.  Then $\mu^{P}$ and $\mu^{Q}$ are the joinings $(\mu^{i_1,\dots, i_{k-2}})^{U',V'}, (\mu^{i_1,\dots, i_{k-2}})^{U',V'},$ and equality follows from the case $k=2,$ so $\mu^{(i_1,\dots, i_k)}$ does not depend on the order of the indices $i_1,\dots, i_k.$   $\square$

Now we can define, unambiguously, $\mu^{\eta}$ for a subset $\eta=\{i_1,\dots,i_k\}$ of $\{1,\dots, d\},$ as the joining $\mu^{(i_1,\dots, i_k)}.$ 

Given $\eta\subset \{1,\dots, k\},$ we can define a seminorm $\varnorm{\cdot}_\eta$ on $L^\infty(X)$ by 
$$
\varnorm{f}_\eta^{2^k}=\int \bigotimes_{\Ge\in \{0,1\}^k} f(x_\Ge)\, d\mu^\eta(\bx).
$$  One checks that this is a seminorm, and that $\int \bigotimes_{\Ge\in \{0,1\}^k} f_\Ge\, d\mu^{\eta}\leq \min_{\Ge\in \{0,1\}^k}\varnorm{f_\Ge}_{\eta},$ just as in the proof of Proposition 2 of \cite{H}.

The next lemma is Lemma 2 of \cite{H}, adapted to the setting of actions of amenable groups.

We will use the main inequality in the proof of the van der Corput trick from \cite{BMZ}, encoded as the following observation.  (See also \cite{BR}.)

\begin{observation}
 Let $\Phi$ and $\Psi$ left F{\o}lner sequences in $G,$ and let $c>0.$  For each $m\in \BN,$ Let $N(m)$ be large enough that $|\Phi_{N(m)}\triangle g\Phi_{N(m)}|<c|\Phi_{N(m)}|$ for all $g\in \Psi_{N(m)}.$  If $\{x_g\}_{g\in G}$ is a sequence of vectors in a Hilbert space $\MH$ with $\|x_g\|\leq 1$ for all $g,$ we then have
$$
\left\|\frac{1}{|\Phi_{N(m)}|}\sum_{g\in \Phi_{N(m)}} x_g- \frac{1}{|\Psi_m|} \sum_{h\in \Psi_m} \frac{1}{|\Phi_{N(m)}|}\sum_{g\in \Phi_{N(m)}} x_{hg}\right\|<2c.
$$
\end{observation}

\begin{lemma}
Let $\eta\subset \{1,\dots, d\}, |\eta|=k.$ Let $\Phi^{(1)},\dots, \Phi^{(k)}$ be F{\o}lner sequences in $G,$ and let $f_{\Ge}\in L^\infty(X), \Ge\in \{0,1\}^k$ be uniformly bounded by $1.$   Then for all $\delta>0,$ there exists $N_0$ such that whenever $N_i>N_0,$
$$
\left|\frac{1}{|\Phi^{(1)}_{N_1}|}\cdots \frac{1}{|\Phi^{(k)}_{N_k}|} \sum_{g_i \in \Phi_i} \int \prod_{\varepsilon \in \{0,1\}^k} f_\Ge \circ  \prod_{i \in \eta} (T_g^{(i)})^{\Ge_i}\, d\mu \right| \leq \varnorm{f_{\bf 0}}_{\eta}+\delta.
$$
\end{lemma}

\ti{Proof.}  Let 
$$
J=\frac{1}{|\Phi^{(1)}_{N_1}|}\cdots \frac{1}{|\Phi^{(k)}_{N_k}|} \sum_{g_i \in \Phi_i} \int \prod_{\varepsilon \in \{0,1\}^k} f_\Ge \circ  \prod_{i \in \eta} (T_g^{(i)})^{\Ge_i}\, d\mu.
$$
Let $\Psi$ be a F{\o}lner sequence in $G,$ let $\Ga>0,$ and for each $m,$ choose $N(m)$ so that $|\Phi_{N}^{(i)}\triangle g\Phi_{N}^{(i)}|,|\Phi_{N}^{(i)}\triangle g\Phi_{N}^{(i)}|<\Ga|\Phi_{N}^{(i)}|,$ for all $i\leq k,$ and all $g\in \Psi_m,$ whenever $N>N(m).$

Fix some $M>0,$ and fix $N_1,\dots, N_k> \max_{m\leq M} N(m).$

Write $\prod_{\varepsilon \in \{0,1\}^k} f_\Ge \circ  \prod_{i \in \eta} (T_{g_i}^{(i)})^{\Ge_i}$ as 
\begin{align}\label{display}
\prod_{\Ge\in \{0,1\}^k, \Ge_k=0} f_\Ge \circ  \prod_{i \in \eta} (T_{g_i}^{(i)})^{\Ge_i} \cdot \prod_{\Ge\in \{0,1\}^k, \Ge_k=1} f_\Ge \circ  \prod_{i \in \eta} (T_{g_i}^{(i)})^{\Ge_i}.
\end{align}
Writing $F_{g_1,\dots,g_{k-1}}$ for the first factor in (\ref{display}), and $F_{g_1,\dots, g_k}'\circ T_{g_k}^{(k)}$ for the second factor,  we have
$$
J= \frac{1}{|\Phi^{(1)}_{N_1}|\cdots |\Phi^{(k-1)}_{N_{k-1}}|} \sum_{g_i\in \Phi_{N_i}^{(i)}} \int \frac{1}{|\Phi_k|} \sum_{g_k\in |\Phi_k|} F_{g_1,\dots, g_{k-1}}\cdot  F'_{g_1,\dots, g_{k-1}}\circ T_{g_k}^{(k)} \, d\mu.
$$
Applying the Cauchy-Schwartz inequality, we find
$$
|J|^2\leq \frac{1}{|\Phi^{(1)}_{N_1}|\cdots |\Phi^{(k-1)}_{N_{k-1}}|} \sum_{g_i\in \Phi_{N_i}^{(i)}} \left|\int \frac{1}{|\Phi_k|} \sum_{g_k\in |\Phi_k|} F'_{g_1,\dots, g_{k-1}}\circ  T_{g_k}^{(k)} \, d\mu\right|^2_{L^2(\mu)}.
$$
By the observation above, replacing the average inside the norm by 
\begin{align}\label{replacement}
\frac{1}{|\Psi_m|}\sum_{h_k\in \Psi_m}\int \frac{1}{|\Phi_k|} \sum_{g_k\in |\Phi_k|}   F'_{g_1,\dots, g_{k-1}}\circ T^{(k)}_{h_kg_k}\, d\mu
\end{align}
introduces an error of at most $2\Ga,$ so we estimate the norm of (\ref{replacement}).    The square of the norm is at most
$$
\frac{1}{|\Psi_m|^2}\sum_{g_k, j_k \in \Psi_m} \int F'_{g_1,\dots,g_{k-1}}\circ T^{(k)}_{h_kg_k} \cdot F'_{g_1,\dots, g_{k-1}}\circ T^{(k)}_{j_kg_k}\, d\mu.
$$
Applying $T^{(k)}_{g_k^{-1}j_k^{-1}}$  to the integrand, this becomes
$$
\frac{1}{|\Psi_m|^2} \sum_{g_k h_k,\in \Psi_{m}} \int F'_{g_1,\dots,g_{k-1}} \circ T^{(k)}_{h_kj_k^{-1}} \cdot F'_{g_1,\dots, g_{k-1}}\, d\mu.
$$

Then $|J|^2$ is bounded by
\begin{align*}
\Ga+\frac{1}{|\Phi_{N_1}^{(1)}|\cdots |\Phi_{N_{k-1}}^{(k-1)}|}\sum_{\substack{g_1\in \Phi_{N_1}^{(1)} \\ \cdots\\ g_{k-1}\in \Phi_{N_{k-1}}^{(k-1)} } } \frac{1}{|\Psi_{m_k}|^2}\sum_{j, h\in \Psi_{m_k}} \int F'_{g_1,\dots,g_{k-1}}\circ T^{(k)}_{hj^{-1}}  \cdot F'_{g_1,\dots, g_{k-1}}\, d\mu.
\end{align*}
Repeating this argument $k$ times produces the inequality
\begin{align}\label{bound}
|J|^{2^k}< C\alpha +\frac{1}{|\Psi_{m_1}|\cdots |\Psi_{m_k}|}\sum_{\substack{j_1, h_1\in \Psi_{m_1}^{(1)}\\ \cdots \\  j_k, h_k\in \Psi_{m_k}^{(k)} }} \int \prod_{\Ge\in \{0,1\}^k} f_{\mathbf 0} \circ \prod_{i=1}^{k} (T^{(i)}_{h_i j_i^{-1}})^{\Ge_i}\, d\mu,
\end{align}
where $C$ depends only on $k.$

Now for a given $\delta>0,$ by Proposition \ref{variation} and the definition of $\varnorm{\cdot}_\eta,$ there exists $M$ so that there are choices $m_1,\dots, m_k<M$ making the average on the right-hand side of (\ref{bound}) differ from $\varnorm{f_{\mathbf 0}}_{\eta}^{2^k}$ by at most $\frac{1}{2}\delta.$   Having chosen $\Ga=\frac{1}{2C}\delta,$ we get $J<\varnorm{f_\mathbf 0}_{\eta}^{2^k}+ \delta.$ $\square$

\subsection{The magic extension and proof of part (1) of Theorem \ref{main}.}

In this section we use the joining $\mu^{(1,\dots,d)}$ to construct an extension of the system $\BX=(X,\MX,\mu,T^{(1)},\dots, T^{(d)})$ with some nice properties.

Let $X^*=X^{\{0,1\}^d}, \MX^*= \MX^{\otimes 2^d}, \mu^*=\mu^{(1,\dots, d)},$ and for $1\leq i \leq d,$ define $T^{(i)*}:G\times X^*\to X^*$ by
$$
T^{(i)*}_gx=\left\{\begin{aligned} T_g^{(i)} x_\Ge &\text{ if } \Ge_i=0;\\
 x_\Ge &\text{ if } \Ge_i=1.  \end{aligned}\right.
$$
One can check from the definition of $\mu^*$ that it is invariant under each $T^{(i)}.$  The system $\mathbf X^*=(X^*,\MX^*,\mu^*, T^{(1)*}, \dots, T^{(d)*})$ is called the \ti{magic} extension of $\BX.$

Furthermore, $\BX$ is a factor of $\BX^*,$ with the factor map defined as projection on the first coordinate. 

Let $\mu^{**}$ be the measure constructed from $\mathbf X^*$  the way $\mu^*$ was constructed from $\BX.$  For $\eta\subset \{1,\dots, d\},$ let $\varnorm{\cdot}_\eta^{*}$ be the seminorm on $L^\infty(\mu^*)$ constructed from  $\mathbf X^*$ the way $\varnorm{\cdot}_\eta$ was constructed from $\mathbf X.$   

For $\eta\subset \{1,\dots, d\},$ let $\MZ_{\eta}$ denote the $\sigma$-algebra spanned by the $T^{(i)}$-invariant sets for $i\in \eta.$  That is,
$$
\MZ_{\eta}:= \bigvee_{i\in \eta} \MI_i.
$$

Exactly as in Theorem 3.2 of \cite{C}, we have
\begin{theorem}
For every $\Ge\subseteq \{1,\dots, d\}, \Ge\neq \emptyset,$ and every function $f\in L^\infty(\mu^*),$ if $\BE_{\mu^*}(f|\MZ_\Ge^*)=0$ then $\varnorm{f}_\Ge^*=0.$
\end{theorem}
The proof is identical to the proof in \cite{C}, so we omit it.

We may proceed as in the proof of Proposition 3.3 in \cite{C} to prove that
\begin{proposition}\label{C3.3}  Let $f_\Ge, \Ge\in \{0,1\}^d$  be functions on $X^{\{0,1\}^d}$ with $\varnorm{f_\Ge}_{L^\infty(\mu^*)} \leq 1$ for every $\Ge.$  Let $\Phi^{(1)},\dots, \Phi^{(d)}$ be F{\o}lner sequences in $G.$  Then the limit
\begin{align}\label{limit*}
\lim_{\min(N_1,\dots,N_d)\to \infty} \frac{1}{|\Phi^{(1)}_{N_1}|\cdots |\Phi^{(d)}_{N_d}|} \sum_{g\in \Phi_{N_1}^{(1)}\times \cdots \times \Phi_{N_d}^{(d)}} \prod_{i=1}^d f_\Ge \circ (T^{(i)*}_{g_i})^{\Ge_i}
\end{align}
exists in $L^2(\mu^*).$
\end{proposition}
Since $\BX$ is a factor of $\BX^*$, we can consider $L^\infty(\mu)$ as a subalgebra of $L^\infty(\mu^*),$ and conclude that the convergence asserted in part (1) of Theorem \ref{main} holds.

\begin{remark} For applications, it is useful to know that the limit is independent of the choice of F{\o}lner sequences $\Phi^{(i)}.$  This can be accomplished by interpolating F{\o}lner sequences: given two collections of F{\o}lner sequences $\Phi^{(i)},\Psi^{(i)},i=1,\dots, d$ let $\Theta^{(i)}_{2N}= \Phi^{(i)}_N, \Theta^{(i)}_{2N+1}=\Psi^{(i)}_N.$  Then the limit in (\ref{limit*}) involving the F{\o}lner sequence $\Psi^{(i)}$ must agree with the limit involving $\Theta^{(i)}$ and the limit involving $\Phi^{(i)}.$
\end{remark}

\subsection{Proof of the inequality in Theorem \ref{main}.}

In Section 4 of \cite{Gr}, it is shown that if $\{u_{g,h}\}_{g,h\in G}$ is a bounded sequence indexed by $G\times G$ such that for all two-sided F{\o}lner sequences $\Phi,\Psi,$ the limits
\begin{align*}
L&=\lim_{N\to\infty} \frac{1}{|\Phi_N||\Psi_N|}\sum_{(g,h)\in \Phi_N\times\Psi_N} u_{g,h}\\
L'&= \lim_{N\to \infty} \frac1{|\Phi_N|}\sum_{g\in \Phi_N} \lim_{M\to \infty} \frac{1}{|\Psi_M|}\sum_{h\in \Psi_M} u_{g,h}
\end{align*}
exist and are independent of the choice of F{\o}lner sequences $\Phi,\Psi,$  then $L=L'.$  By an identical argument, we have that given a sequence $\{u_g\}_{g\in G^d},$ indexed by $G^d,$ if the limits
\begin{align*}
J&=\lim_{N\to \infty} \frac{1}{|\Phi^{(1)}_N|\cdots |\Phi^{(d)}_N|}\sum_{g\in \Phi^{(1)}_N\times \cdots \times \Phi^{(d)}_N} u_{g}\\
J'&=\lim_{N_1\to \infty} \frac{1}{|\Phi^{(1)}_N|} \sum_{g_1\in \Phi^{(1)}_{N_1}} \cdots \lim_{N_d\to \infty} \frac{1}{|\Phi^{(d)}_N|} \sum_{g_d\in \Phi^{(d)}_{N_d}} u_{(g_1,\dots, g_d)}
\end{align*}
exist and are independent of the choice of F{\o}lner sequences involved, then $J=J'.$  

\ti{Proof of Theorem \ref{main}, part (2).}

The inequality in Theorem \ref{main} is proved inductively in the case $d=2$ in \cite{Gr}, and the induction may be continued to achieve a proof for arbitrary $d.$  

Suppose that the inequality holds for $d-1$ rather than $d.$  
Writing $F_{g_1,\dots, g_{d-1}}$ for $\prod_{\Ge\in \{0,1\}^d, \Ge_d=0} f\circ \prod_{i=1}^d(T^{(i)}_{g_i})^{\Ge_i},$ we have
$$
\prod_{\Ge\in \{0,1\}^d} f\circ \prod_{i=1}^d (T^{(i)}_{g_i})^{\Ge_i}= F_{g_1,\dots, g_{d-1}} \cdot F_{g_1,\dots, g_{d-1}}\circ T_{g_d}^{(d)}.
$$
Integrating and averaging over $\Phi^{(d)}$, we have
\begin{align*}
\lim_{N\to \infty} \frac{1}{|\Phi^{(d)}_N|} \sum_{g\in \Phi^{(d)}_N} \int F_{g_1,\dots, g_{d-1}} \cdot F_{g_1,\dots, g_{d-1}}\circ T_{g_d}^{(d)}& \, d\mu \\ &= \int  \BE(F_{g_1,\dots, g_{d-1}}|\MI_d)^2  \, d\mu\\
&\geq \left(\int\BE(F_{g_1,\dots, g_{d-1}}|\MI_d)\, d\mu  \right)^2\\
&= \left(\int F_{g_1,\dots, g_{d-1}}\, d\mu  \right)^2.
\end{align*}
Averaging over $\Phi^{(1)}\times \cdots \times \Phi^{(d-1)},$ and applying the induction hypothesis, we find
\begin{align*}
\frac{1}{|\Phi^{(1)}_N|\cdots |\Phi^{(d-1)}_N|} \sum_{g_i\in \Phi^{(i)}_N}& \left(\int  F_{g_1,\dots, g_{d-1}}\, d\mu  \right)^2   \geq \\ & \left(\frac{1}{|\Phi^{(1)}_N|\cdots |\Phi^{(d-1)}_N|}\sum_{g_i\in \Phi^{(i)}_N}\int F_{g_1,\dots, g_{d-1}}\, d\mu  \right)^2.
\end{align*}
Averaging and applying the induction hypothesis yields the desired inequality.  $\square$

\begin{observation} Part (2) of Theorem \ref{main} and Lemma \ref{tosynd} imply that for all $c>0$
$$
R:=\{(g_1,\dots,g_d)\in G^d: \int \prod_{\Ge\in \{0,1\}^d} f\circ \prod_{i=1}^d (T^{(i)}_{g_i})^{\Ge_i}\, d\mu> \left(\int f \, d\mu \right)^{2^d}-c\}
$$
is both left- and right syndetic, as $R$ must meet  at least one element of each sequence of the form $\{\Phi^{(1)}_N\times \cdots \times \Phi^{(d)}_N\}_{N\in \BN},$ where each $\Phi^{(i)}$ is a F{\o}lner sequence in $G.$
\end{observation}

To prove Corollary \ref{maincor}, we apply a variation of the Furstenberg correspondence principle as stated in \cite{BMZ}.  Let $G$ be a countable amenable and consider the shift space $\{0,1\}^{G^d}.$ Let $T^{(i)}$ denote the left shift in $i$th coordinate of $G^d,$ so that $(T^{(i)}_g\xi)(g_1,\dots, g_i,\dots, g_d)=\xi(g_1,\dots, g^{-1}g_i,\dots, g_d).$

\begin{prop}\label{correspondence}  Let $\Psi$ be a left F{\o}lner sequence in $G^d.$  Suppose $E\subseteq G^{d}$ with $d_\Psi(E)>0.$    Let $X=\overline{\{\prod_{i=1}^{d} T_{g_i}^{(i)}1_E: g_i\in G\}}$ be the orbit closure of $1_E$ in $\{0,1\}^{G^d}.$  If
\begin{align*}
d_{\Psi}(E)=\limsup_{n\to \infty} \frac{|E\cap \Psi_N|}{|\Psi_N|}>0,
\end{align*}
then there exists a probability measure $\mu$ on $X,$ invariant under each $T^{(i)},$ such that
\begin{align*}
\mu(\{\xi\in X:\xi(e,\dots, e)=1\})\geq d_\Psi(E).
\end{align*}
Let $A:=\{\xi\in X:\xi(e,\dots, e)=1\}.$ For all collections $\{g_{i,j}: 1\leq i \leq d, 1\leq j \leq n\}$ of elements of $G,$ the inequality
$$
d_\Psi\left(\bigcap_{j=1}^n (g_{1,j},\dots, g_{d,j})E\right )\geq \mu\left(\bigcap_{j=1}^n \prod_{i=1}^d T^{(i)}_{g_{i,j}} A\right)
$$
holds.
\end{prop}
(This is stated in \cite{BMZ} for the case $d=2,$ but a similar proof works for arbitrary $d.$  Furthermore, the inequality $\mu(\{\xi\in X:\xi(e,\dots, e)=1\})>0$ is stated in \cite{BMZ}, but the inequality $\mu(\{\xi\in X:\xi(e,\dots, e)=1\})\geq d_\Psi(E)$ is actually proved.  Also, the right shift, rather than the left shift, is used in \cite{BMZ}.)

\ti{Proof of Corollary \ref{maincor}.}  Let $\Psi$ be a left F{\o}lner sequence in $G^d,$ and let $E\subset G$ with $d_\Psi(E)=\delta>0.$  Let $T^{(i)},i=1,\dots, d,$ $A$ and $\mu$ be as in Proposition \ref{correspondence}, so that $\int 1_A\, d\mu=\delta.$  

Let $c>0.$  By Theorem \ref{main}, the set 
$$\left\{(g_1,\dots, g_d): \int \prod_{\Ge\in \{0,1\}^d} 1_A\circ \prod_{i=1}^d(T^{(i)}_{g_i})^{\Ge_i}\, d\mu> \delta^{2^d}-c\right\}
$$
is both left- and right syndetic in $G^d.$  By Proposition \ref{correspondence}, this implies that the set of $(g_1,\dots, g_d)$ such that $d_\Phi(\bigcap_{\Ge\in\{0,1\}^d} (g_1^{-\Ge_1},\dots, g_d^{-\Ge_d})E)>\delta^{2^d}-c$ is both left- and right syndetic in $G^d.$ 
$\square$

\end{document}